\documentclass[11pt]{amsart}

\usepackage{amsmath,amssymb,fullpage,color,graphicx}

\numberwithin{equation}{section}

\theoremstyle{plain}

\newtheorem{theorem}[equation]{Theorem}
\newtheorem{conj}[equation]{Conjecture}

\theoremstyle{remark}

\newtheorem{remark}[equation]{Remark}

\theoremstyle{definition}

\newcommand{\C}{\mathbb C}

\newcommand{\PP}{\mathbb P}
\newcommand{\Q}{\mathbb Q}
\newcommand{\R}{\mathbb R}
\newcommand{\T}{\mathbb T}
\newcommand{\Z}{\mathbb Z}

\newcommand{\fraka}{\mathfrak a}

\newcommand{\frakd}{\mathfrak d}
\newcommand{\frakD}{\mathfrak D}

\newcommand{\frakn}{\mathfrak n}
\newcommand{\frakM}{\mathfrak M}
\newcommand{\frakN}{\mathfrak N}
\newcommand{\frakp}{\mathfrak p}

\newcommand{\calO}{\mathcal O}

\DeclareMathOperator{\Coind}{Coind}
\DeclareMathOperator{\End}{End}
\DeclareMathOperator{\Gal}{Gal}
\DeclareMathOperator{\GL}{GL}

\DeclareMathOperator{\M}{M}
\DeclareMathOperator{\nrd}{nrd}
\DeclareMathOperator{\PSL}{PSL}
\DeclareMathOperator{\Tr}{Tr}
\DeclareMathOperator{\Nm}{N}

\newcommand{\defi}{\textsf}

\begin{document}

\title[A database of Hilbert modular forms]{A database of Hilbert modular forms}
\author{Steve Donnelly}
\address{School of Mathematics and Statistics F07 \\ University of Sydney \\ Sydney, NSW\ 2006,  Australia} \email{donnelly@maths.usyd.edu.au}
\author{John Voight}
\address{Department of Mathematics, Dartmouth College, 6188 Kemeny Hall, Hanover, NH 03755, USA}
\email{jvoight@gmail.com}

\begin{abstract}
We describe the computation of tables of Hilbert modular forms of parallel weight $2$ over totally real fields.
\end{abstract}

\maketitle

Tables of classical modular forms have been computed for many decades now, starting with a table of rational eigenforms (corresponding to modular elliptic curves) appearing in the proceedings of the Antwerp conference in 1972 \cite{AntwerpIV}.  Large tables of classical modular forms now exist, as computed by Cremona \cite{Cremona0,CremonaTables} for rational eigenforms, Stein \cite{Stein,SteinModforms}, and more recently appearing in the \emph{$L$-functions and Modular Forms Database} (LMFDB) \cite{LMFDB}, computed by Stephan Ehlen and Fredrik Str\"omberg.  In this article, we continue this tradition of tabulating automorphic forms on Shimura varieties by computing extensive tables of Hilbert modular cusp forms of parallel weight $2$ over totally real fields.  This work is in a sense complementary to more recent work \cite{Steinetal} on tables of (modular) elliptic curves over $\Q(\sqrt{5})$.  Our data is available at the second author's webpage \cite{Voight-hmf} and is searchable online at the LMFDB \cite{LMFDB}.

In Section 1, we give an overview of the basic algorithms employed.  In Section 2, we describe in more detail the implementation in the computer algebra system \textsc{Magma} \cite{Magma}, providing some documentation and detail about computations in practice and a few auxiliary algorithms employed.  In Section 3, we gather the data which was computed and discuss the extent of the computation.  Finally, in Section 4 we some interesting aspects of the data in detail.

The authors would like to thank Lassina Demb\'el\'e and Matthew Greenberg for useful discussions.  Thanks also go to John Cremona and Aurel Page for adding Galois conjugate Hilbert modular form data to the LMFDB.  The second author was supported by an NSF Grant (DMS-0901971) during the time that these computations were first undertaken, and by an NSF CAREER Award (DMS-1151047) while work was completed.

\section{The algorithm: an overview}

In this section, we give a basic overview of the algorithm we employ.  For a motivated introduction, see the expository work by Demb\'el\'e--Voight \cite{DembVoight}; for a more basic reference on Hilbert modular forms, see e.g.\ Freitag \cite{Freitag}.

Throughout, let $F$ be a totally real field of degree $n=[F:\Q]$ with ring of integers $\Z_F$ and let $\frakN \subseteq \Z_F$ be an ideal.  Let $S_k(\frakN)$ denote the space of Hilbert cusp forms of level $\frakN$ and parallel weight $k \in \Z_{\geq 2}$ for $F$, and let $S_k^{\textup{new}}(\frakN)$ be the new subspace.

The algorithm we use combines methods of Demb\'el\'e \cite{Dembele} and Demb\'el\'e--Donnelly \cite{DD}, along with methods of Greenberg--Voight \cite{GV} and Voight \cite{VoightANTS}.

\begin{theorem} \label{algthm}
There exists an (explicit) algorithm that, on input a totally real field $F$ and an ideal $\frakN\subseteq \Z_F$, computes as output the systems of eigenvalues for the Hecke operators $T_\frakp$ with $\frakp\nmid \frakN$ and the Atkin-Lehner involutions $W_{\frakp^e}$ with $\frakp^e\parallel \frakN$ that occur in the space $S_k^{\textup{new}}(\frakN)$.
\end{theorem}

In other words, there exists an explicit finite procedure which takes as input the totally real field $F$, an ideal $\frakN \subseteq \Z_F$ (encoded in bits in the usual way \cite{Cohen}), and an even integer $k \geq 2$, and outputs a finite set of number fields $E_f \subseteq \overline{\Q}$ and sequences $(a_f(\frakp))_\frakp$ encoding the Hecke eigenvalues for each eigenform constituent $f$ in $S_k(\frakN)$ and $a_f(\frakp) \in E_f$. 

In this article, we consider the first and essential case $k=2$; we believe it would be an interesting project to extend these efforts to include tabulation of Hilbert modular forms of higher weight as well.

The general strategy of the algorithm in Theorem \ref{algthm} is as follows.  Let $B$ be the quaternion algebra over $F$ which is unramified at all finite places and ramified at either all or all but one real place according as $n=[F:\Q]$ is even or odd.  Then by the Eichler--Shimizu--Jacquet--Langlands correspondence, there exists an isomorphism of Hecke modules
\[ S_2^B(\frakN) \xrightarrow{\sim} S_2(\frakN) \]
where $S_2^B(\frakN)$ denotes the space of quaternionic modular forms over $B$ of weight $k$ and level $\frakN$.  The description of the space $S_2^B(\frakN)$ varies accordingly as $n$ is even or odd.  Let $\calO \subseteq B$ be a maximal order.

If $n$ is even, we use Brandt matrices, following the method of Demb\'el\'e \cite{Dembele} and Demb\'ele--Donnelly \cite{DD}.  The Hecke module $S_2^B(\frakN)$ is computed via an explicit action on maps from the set of right ideal classes in $\calO$ to $\PP^1(\Z_F/\frakN)$.  We refer to this as the \defi{definite method}, since the algebra $B$ in this case is definite.

If instead $n$ is odd, the algorithm in Theorem \ref{algthm} is supplied by Greenberg--Voight \cite{GV} and Voight \cite{VoightANTS}.  For simplicity of exposition, assume $F$ has strict class number $1$; the method works for arbitrary class number.  Let $\calO_1^\times=\{\gamma \in \calO : \nrd(\gamma)=1\}$ denote the group of units of $\calO$ with reduced norm $1$.  Suppose that the unique real split place of $B$ corresponds to the embedding $\iota_\infty : B \hookrightarrow \M_2(\R)$ and let $\Gamma(1) = \iota_\infty( \calO_1^\times)/\{ \pm 1 \} \subseteq \PSL_2(\R)$.  Let $\calO_0(\frakN)$ be an Eichler order of level $\frakN$ and let $\Gamma_0(\frakN) = \iota_\infty( \calO_0(\frakN)_1^\times)/\{\pm 1\}$.  Then by the theorem of Eichler-Shimura, and Shapiro's lemma, we have a further isomorphism of Hecke modules
\[ S_2^B(\frakN) \xrightarrow{\sim} H^1\bigl(\Gamma(1),\Coind_{\Gamma_0(\frakN)}^{\Gamma(1)} \C\bigr)^+ \]
where $+$ denotes the $+1$-eigenspace for complex conjugation.  The Hecke module $S_2^B(\frakN)$ is then computed via this isomorphism as group cohomology, using an explicit presentation for the group $\Gamma(1)$ based on the computation of a fundamental domain \cite{Voight-funddom}.  We call this the 
\defi{indefinite method}.  

\begin{remark} \label{JLrmk}
More generally, one can use a quaternion algebra $B$ over $F$ of discriminant $\frakD$, ramified at all or all but one real place, and $\frakM$ an ideal coprime to $\frakD$; then with $\frakN=\frakD\frakM$, the Jacquet-Langlands correspondence yields an isomorphism
\[ S_2^B(\frakM) \xrightarrow{\sim} S_2(\frakN)^{\frakD\textup{-new}} \]
where $S_2(\frakN)^{\frakD\textup{-new}} \subseteq S_2(\frakN)$ denotes the space of Hilbert cusp forms which are new at all primes dividing $\frakD$.  Since $\frakD$ is squarefree, these methods overlap when there is a prime $\frakp$ which exactly divides the level $\frakN$.  We use this to verify our calculations by computing with both the definite and indefinite method when it applies; for more, see Section 3.
\end{remark}

\section{Computing systems of Hecke eigenvalues}

In this section, we discuss the \textsc{Magma} \cite{Magma} implementation of the above algorithms.  These methods saw a first implementation over $\Q$ by David Kohel, and they rely heavily on code for quaternion algebras over number fields, due to Donnelly, Markus Kirschmer, Nicole Sutherland, and Voight (documented in work of Kirschmer--Voight \cite{KV}).

Our algorithm has 6 main steps:
\begin{enumerate}
\item[1.] Precomputation
\item[2.] Computing a basis
\item[3.] Hecke computation via reduction
\item[4.] Decomposition into newforms
\end{enumerate}

We refer to Demb\'el\'e--Voight \cite{DembVoight} for an explicit mathematical description, and here we provide some additional comments about the implementation.  

\subsection*{Step 1: Precomputation}

The first step of the algorithm involves precomputation steps that only depend on the field $F$: that is to say, for all levels $\frakN$ with $\frakp \parallel \frakN$, this step need only be performed once.  For the definite method, the precomputation step involves determining a set of representatives for the right ideal classes in $\calO$ up to isomorphism.  

\begin{itemize}
\item We precompute the theta series of each left order $\calO$ in the set of ideal representatives, recording 
\[ \#\{\alpha \in \calO_L(I) : \Tr_{F/\Q} \nrd(\alpha) = n\} \] 
for a small set of values $n \in \Z_{\geq 0}$.  This reduces the number of isomorphism tests of ideals when computing Hecke operators below, but at the marginal additional cost of requiring the computation of left orders and representation numbers.  
\item As computing inverses and colon ideals can also be time consuming, we compute for each representative right ideal $I$ a small element $b \in I$ and a left ideal $I'$ of prime power norm such that $\calO_L(I) b = II'$: then $J I' = J (I^{-1})b$, and $\alpha$ generates $JI'$ if and only if $\alpha/b$ generates $JI^{-1}$.
\end{itemize}

For the indefinite method, we compute a fundamental domain for the Shimura curve associated to the maximal order, together with a presentation for the group $\Gamma_0(1)$ \cite{Voight-funddom}.

\subsection*{Step 2: Computing a basis}

To compute a basis of forms, we need fast algorithms for working with $\PP^1(\Z_F/\frakN)$, and our efficient implementation treats the case where $\frakN$ is prime and composite differently: in the latter we index elements $[a:b] \in \PP^1(\Z_F/\frakN)$ according to divisor $\fraka=(a) \mid \frakN$.  With these representatives in hand, in the definite case we list orbits of $\PP^1(\Z_F/\frakN)$ under the action of the stabilizer groups and in the indefinite case we record the permutation action in the coinduced module.  (If we were working in higher weight, and a nontrivial representation $V$, we would put a copy of $V$ on each element of $\PP^1(\Z_F/\frakN)$.)  

For sanity, we also verify the dimension of the space using a dimension formula (involving class numbers and embedding numbers).

\subsection*{Step 3: Hecke computation via reduction}

For the computation of Hecke operators, we use reduction on $\PP^1(\Z_F/\frakN)$ in the definite case and reduction theory in the fundamental domain in the indefinite case.

\subsection*{Step 4: Decomposition into newforms}

With Hecke matrices in hand, we now outline the procedures used in our implementation to decompose a newspace as a direct sum of new simultaneous eigenspaces.  Let $S=S_2(\frakN)$.  We use the following properties.

\begin{enumerate}
\item The Hecke operators $T_\frakp$ for $\frakp$ coprime to $\frakN$ all commute.
\item The Hecke eigenspaces in $S$ have distinct systems of eigenvalues for the operators $T_\frakp$ with $\frakp$ coprime to $\frakN$ (``multiplicity one'').
\item We have $\dim S^{\textup{$\frakp$-old}} = 2\dim S_2(\frakN/\frakp) - \dim S_2(\frakN/\frakp^2)$, with the last term omitted if $\frakp^2$ does not divide $\frakN$.
\end{enumerate}

These properties are standard: for (2), see for instance work of Shemanske and Walling \cite[Theorem 3.6]{ShemanskeWalling}, and (3) follows from standard facts about degeneracy maps (the proof is a straightforward generalization of the classical proof (where $F=\Q$), and is omitted).

Using (1) and standard linear algebra techniques, we decompose the space into Hecke constituents corresponding to $\Q$-irreducible subspaces (a direct sum of eigenspaces).  As harmless as it may seem, this step is often the bottleneck in computation.

Finally, we show how to decompose a space into old and new subspaces.  For a prime $\frakp \nmid \frakN$, the $\frakp$-old subspace $S^{\textup{$\frakp$-old}}$ of $S=S_2(\frakN)$ is the span of the images of all degeneracy maps $S_2(\frakN/\frakp) \to S_2(\frakN)$.  The old subspace $S^{\textup{old}}$ is the sum of the spaces $S^{\textup{$\frakp$-old}}$ for primes $\frakp \mid \frakN$.  The newspaces $S^{\textup{new}}$ and $S^{\textup{$\frakp$-new}}$ are the unique Hecke-invariant complements of $S^{\textup{old}}$ and $S^{\textup{$\frakp$-old}}$, respectively.  

Now suppose we are given two Hecke-modules $S$ and $S_0$ (e.g.\ $S_0=S_2(\frakN/\frakp)$ and we wish to compute the Hecke-submodule $S^{\textup{old}}$ of $S$ consisting of all the irreducible Hecke-submodules (i.e.\ simultaneous eigenforms) of $S$ that arise from $S_0$ via degeneracy maps.  In addition, suppose we know a priori the dimension of $S^{\textup{old}}$, as in Property (3).  For a prime $\frakp$ coprime to $\frakN$, the space $S^{\textup{old}}$ is contained in the kernel of the characteristic polynomial of $T_\frakp \mid S_0$ evaluated at $T_\frakp \mid S$.  Since we know the dimension, we may compute $S^{\textup{old}}$ in this manner using finitely many primes $\frakp$, stopping when the intersection has the right dimension.  One could obtain the complement $S^{\textup{new}}$ by computing the Hecke-invariant complement, requiring extra module-splitting effort.  However, for present purposes, we do not need to know $S^{\textup{new}}$ as a particular subspace of $S$; we merely need the Hecke action on $S^{new}$, which can be obtained using the quotient $S/S^{\textup{old}} \cong S^{\textup{new}}$.

Alternatively, we can compute the degeneracy operators themselves, and then compute the orthogonal complement.  

Now suppose that $S$ is a new subspace of dimension $d$.  By (1), it contains $d$ simultaneous eigenvectors having distinct systems of eigenvalues, which we wish to determine.
The Hecke algebra $\T=\Z[T_\frakn]_{\frakn}$ is a commutative $\Z$-algebra which is generated by Hecke operators $T_\frakp$ with $\frakp$ coprime to $\frakN$.  Therefore $\T \otimes \Q$ is a direct product of number fields, an \'etale $\Q$-algebra generated by a single primitive element.  For practical purposes, one may find a primitive element by trial and error (for instance taking simple combinations of the first few Hecke operators), since the set of non-primitive elements must be sparse.  Then the simultaneous eigenspaces of $S$ are simply the eigenvectors associated to the primitive element.

\section{Data computed}

We compute for the totally real fields $F$ of degree $n \leq 6$ with discriminant $\leq 2\cdot 10^n$.  Note that if $[F:\Q] > 6$ then $d_F > 2 \cdot 10^n$ by the Odlyzko bounds; indeed, the totally real field of degree $7$ of smallest discriminant is blank.  

In this range, the discriminant and degree of the field uniquely characterizes the field, with one exceptional pair: in degree $4$, there are two fields of discriminant $16448$: one, with minimal polynomial $x^4-2x^3-6x^2+2$ has Galois group $S_4$ and the other $x^4-2x^3-7x^2+8x+14$ has Galois group $D_8$.  (The first multiple discriminant for cubic fields is $3969$.)

This gives a total of $233676$ Hilbert modular forms in the database.  To save space, we did not include data for conjugate levels; for consistency across the database, these were later added to the LMFDB, increasing the overall size of the data.

The computation was performed using the Vermont Advanced Computing cluster (VACC).  The odd degree fields, which used the indefinite method, took approximately 4 CPU years to compute.  Much of this time is spent computing Atkin-Lehner involutions and coping with higher class number issues.  

The data is hosted at the LMFDB \cite{LMFDB} and backed up at the NECC File Exchange.  The raw data takes 172 GB, compressed 78 GB.

\begin{remark}
As in Remark \ref{JLrmk}, in many situations we can employ either the definite or indefinite method.  For example, if $F$ has odd degree and $\frakp \parallel \frakN$, then one can apply the definite method with a quaternion algebra $B$ ramified at all real places and the prime $\frakp$; the Hecke module computed is then $S_2(\frakN)^{\textup{$\frakp$-new}}$.  To compute a single space of forms, the definite method runs more quickly than the indefinite method, but for tabulation one must use the indefinite method if one wants to compute with forms of square level; it is also more efficient to use the indefinite method since for each quaternion algebra used one must repeat the precomputation step.  
\end{remark}

\begin{center}
\begin{tabular}{ccc|ccc}
$d_F$ & $N$ & \textup{num} & $d_F$ & $N$ & \textup{num} \\
\hline
5 & 4999 & 7583 &     49 & 2059 & 838      \\  
8 & 5000 & 12875 &     81 & 719 & 343         \\
12 & 5000 & 10565 &    148 & 499 & 1068       \\
13 & 1999 & 5837 &    169 & 625 & 515        \\
17 & 988 & 3387 &     229 & 256 & 673       \\
& \vdots &  &   &  \vdots \\
493 & 4 & 28 &      1944 & 54 & 568 \\
497 & 7 & 19 &      1957 & 27 & 199 \\
\hline
& & 104593 & & & 40330 
\end{tabular}

\vspace{3ex}

\begin{tabular}{ccc|ccc|ccc}
$d_F$ & $N$ & \textup{num} & $d_F$ & $N$ & \textup{num} & $d_F$ & $N$ & \textup{num} \\
\hline
725 & 4091 & 2807 &         14641 & 1013 & 217 &       300125 & 911 & 177 \\
1125 & 991 & 375 &         24217 & 835 & 1062 &        371293 & 961 & 251 \\
& \vdots & &         &  \vdots & &           & \vdots & \\
19821 & 97 & 241 &         195829 & 83 & 234 &        1997632 & 127 & 84 \\
\hline
& & 64329 & & & 15402 & & & 9022
\end{tabular}
\end{center}

\section{Comments on data}

\subsubsection*{Base change forms}

We detect base change forms as follows.  First, a necessary condition for $f$ to be a base change from a smaller field $E$ is that $\sigma(\frakN)=\frakN$ and $a_{\sigma(\frakp)}=a_\frakp$ for all $\sigma \in \Gal(F/E)$.  If the form fails this test (for all possible subfields $E$), it is not a base change form; if it passes this test, then is likely to be a base change.  (One could detect more generally forms that are twists of base change forms.)  To be certain, we look in the database for forms over the field $E$ with level supported at primes dividing the discriminant $\frakd_{F/E}$ and $\frakn = \Z_E \cap \frakN$; the precise recipe is provided by Loeffler and Weinstein for forms that are base change from $\Q$.  

\subsubsection*{CM forms}

We detect CM forms as follows.  Let $\frakp$ be a prime which is trivial in the narrow class group; then the Hecke eigenvalue $a_\frakp$ is a totally real algebraic integer and the characteristic polynomial of Frobenius $T^2-a_\frakp T+N(\frakp)$ has discriminant $a_\frakp^2-4N(\frakp)$.  If the form is CM, then the Frobenius belongs to the quadratic field of this discriminant.  If we find two such primes $\frakp$ with $a_\frakp \neq 0$ whose CM discriminants belong to different square classes in $F^\times/F^{\times 2}$, then the form does not have CM; if these square classes agree, then it gives strong likelihood that the form is CM with the given CM discriminant $K$.  One can then also verify that $a_\frakp=0$ for primes $\frakp$ inert in $K$ to provide further evidence.  To conclude that the form really is CM, we could use the machinery for Hecke characters.

\subsubsection*{Computing the $L$-function} 

Given the Hecke eigenvalues $a_\frakp(f)$ of a Hilbert modular form $f$ with trivial central character, we define its $L$-series via the Euler product
\[ L(s,f) = \prod_{\frakp \mid \frakN} (1-a_\frakp (\Nm\!\frakp)^{-s})^{-1} \prod_{\frakp \nmid \frakN} (1-a_\frakp (\Nm\!\frakp)^{-s} + (\Nm\!\frakp)^{-2s})^{-1}. \]
The function $L(s,f)$ has analytic continuation to all of $\C$ and a functional equation
\[ \Lambda(s,f) = (d_F^2 \Nm(\frakN))^{s/2} \Gamma_\C(s+1/2)^{n}L(s,f) = \Lambda(1-s,f). \]

\subsubsection*{Factorization of new subspaces}

The largest Hecke-irreducible space in our database has dimension $286$, and occurs for the field $F=\Q(\sqrt{296})$ and level $\frakN$ of norm $29$.  It would be interesting to consider the distribution of the dimensions of Hecke irreducible subspaces inside the new subspace.




\subsubsection*{Endomorphism algebras}

We can also see some interesting patterns in the distribution of endomorphism algebras of modular abelian surfaces over totally real fields.  First, we have the following conjecture of Coleman (see e.g.\ Bruin--Flynn--Gonz\'alez--Rotger \cite[Conjecture $C(e,g)$]{BFGR}.

\begin{conj} \label{conj:finiteFE}
Let $F$ be a totally real field and let $g \in \Z_{\geq 1}$.  Then the set
\[ \{E=\End(X)_{\Q} : [E:\Q]=g \text{ and $X$ an abelian variety over $F$ with $\dim(X)=g$}\} \]
is finite.
\end{conj}

In particular, there should be only finitely many real quadratic fields $E$ that occur as the endomorphism algebra of a $\GL_2$-type abelian surface ($g=2$) over a fixed totally real field $F$.  (Our data becomes a bit too sparse to say anything for $g \geq 3$ except for a few fields.)  There are a total of $37444$ Hilbert modular forms whose Hecke eigenvalue field has degree $g=2$; corresponding to such a form $f$ over $F$ is an isogeny class of abelian surfaces $X$ over $F$ whose common endomorphism algebra $E=\End(X)_{\Q}$ is a quadratic field.  

Hilbert modular forms over $F$ whose Hecke field $E$ is an imaginary quadratic field correspond to Hecke Grossencharacters.  The largest absolute discriminant of an imaginary quadratic field observed was $E=\Q(\sqrt{-288})$, occurring for a Hilbert modular form over $F=\Q(\sqrt{440})$, a field with class number $2$.  

More interesting are the cases when the Hecke field $E$ is a real quadratic field.  The largest such real quadratic field observed as a Hecke field was $E=\Q(\sqrt{1260})$, occurring over the totally real quartic field $F$ of discriminant $11025$.  The distribution of the quadratic Hecke fields for a fixed field $F$ shows rapid decay by discriminant; the following histogram of the quadratic Hecke fields for $F=\Q(\sqrt{5})$ is typical.

\begin{center}
\includegraphics[scale=0.4]{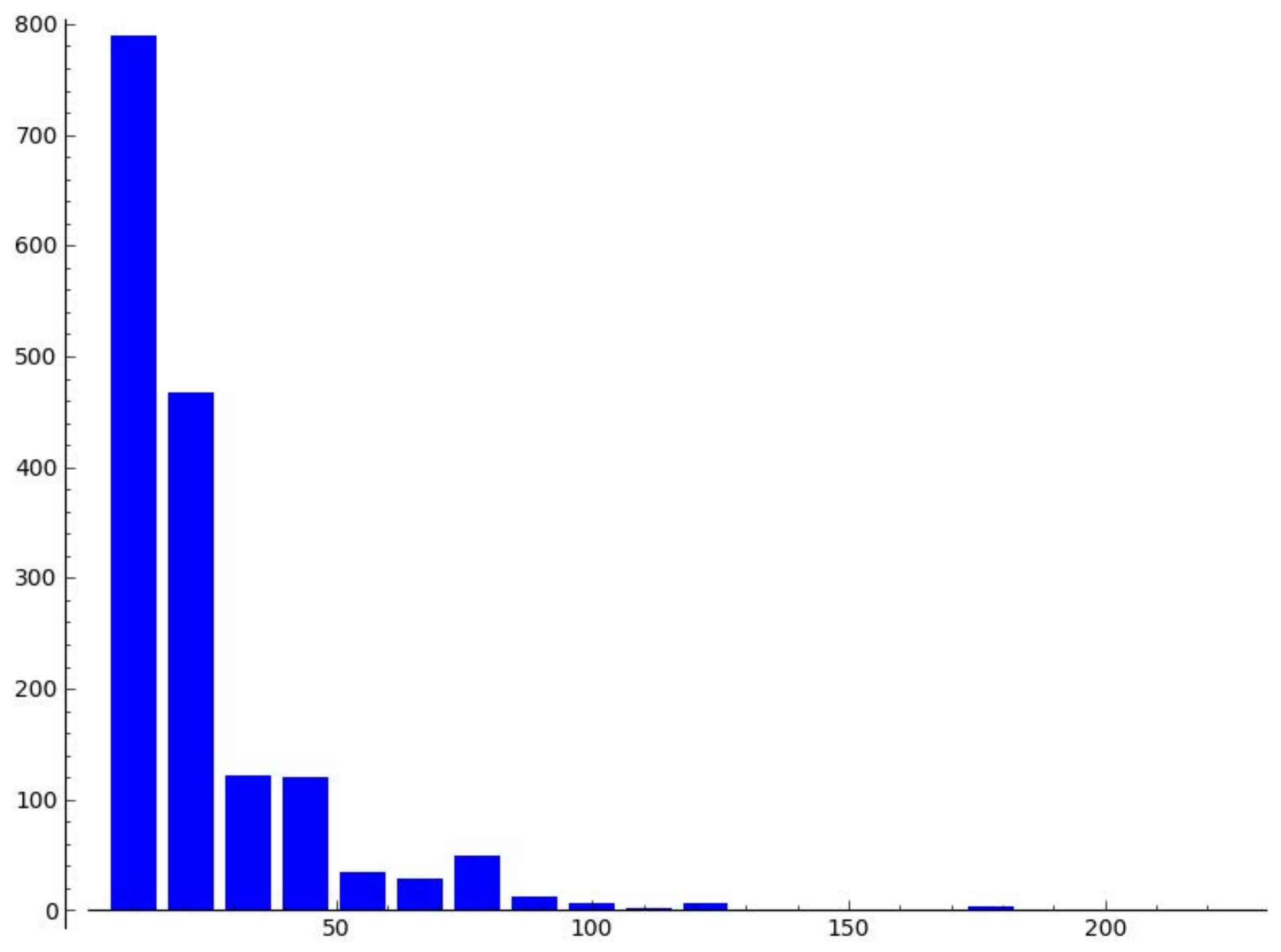}
\end{center}

The largest discriminant field for $F=\Q(\sqrt{5})$ is $E=\Q(\sqrt{228})$, and it occurs just once with level norm $3420$ (up to the computed level norm $4999$).  We tabulate in a similar way the following maximum discriminant quadratic Hecke fields for the first few $F$:
\begin{center}
\begin{tabular}{c|c}
$d_F$ & $\max(d_E)$ \\ \hline
5 & 228 \\
8 & 384 \\
12 & 192 \\
13 & 252 \\
17 & 192 \\
21 & 588 \\
24 & 288 \\
28 & 448
\end{tabular}
\end{center}
Put together, this data gives some weak evidence for Conjecture \ref{conj:finiteFE}---at the very least, it shows that the large discriminant fields occur with very low density.

\end{document}